\documentclass[10pt]{article}
\usepackage{graphicx}
\usepackage{authblk}
\usepackage{mathtools}
\usepackage{amssymb,amsthm}
\usepackage{url}
\usepackage{hyperref}
\usepackage[labelfont={bf}]{caption}
\usepackage{xcolor}

\newtheorem{theorem}{Theorem}

\date{}

\begin{document}
\title{Simultaneous space-time finite element methods for parabolic optimal control problems
   \thanks{Supported by the Austrian Science Fund under the grant W1214, project DK4.}
}

%
\author[1]{Ulrich Langer}
\affil[1]{%
   Institute for Computational Mathematics\authorcr
   Johannes Kepler University Linz\authorcr
   Altenbergerstr. 69, 4040 Linz, Austria}
\author[2]{Andreas Schafelner}
\affil[2]{%
   Doctoral Program ``Computational Mathematics''\authorcr
   Johannes Kepler University Linz\authorcr
   Altenbergerstr. 69, A-4040 Linz, Austria}

\maketitle
\begin{abstract}
   This work presents, analyzes and tests stabilized space-time finite element methods on fully
   unstructured simplicial space-time meshes for the numerical solution of 
   space-time tracking parabolic optimal control problems with the standard
   $L_2$-regularization.\\[1em]
   Keywords: Parabolic optimal control problems; $L_2$-regularization; Space-time finite element methods.
\end{abstract}
%
%
\section{Introduction}
\label{LS:sec:Introduction}
Let us consider the following space-time tracking optimal control problem:
For a given target function $y_d \in L_2(Q)$ (desired state) 
and for some appropriately chosen regularization parameter $\varrho>0$,
find the state 
$y \in Y_0 = \{ v \in L^2(0,T; H_0^1(\Omega)) : \partial_t v \in
   L^2(0,T;H^{-1}(\Omega)), \, v = 0 \mbox{ on }  \Sigma_0\}$
and the control $u \in U = L_2(0,T;L_2(\Omega)) = L_2(Q)$ minimizing the cost functional
\begin{equation}
   \label{LS:eqn:CostFunctional}
   J(y,u) = \frac{1}{2}\int_Q\!|y-y_d|^2\;\mathrm{d}Q + \frac{\varrho}{2} \|u\|_{L_2(Q)}^2 
\end{equation}
subject to 
the linear parabolic initial-boundary value problem (IBVP) 
\begin{equation}
   \label{LS:eqn:IBVP}
   \partial_t y - \mbox{div}_x(\nu \nabla_x y) = u  \mbox{ in } Q,\quad
   y  =  0 \mbox{ on } \Sigma,\quad
   y = 0  \mbox{ on } \Sigma_0,
\end{equation}
where $Q := \Omega \times (0,T)$, $\Sigma:= \partial \Omega \times (0,T)$,
$\Sigma_0:=\Omega\times\{0\}$,  $T > 0$ is the final time, 
$\partial_t$ denotes the partial time derivative, 
$\mbox{div}_x$ is the spatial divergence operator, $\nabla_x$ is the spatial gradient,
and the source term $u$ on the right-hand side of the parabolic PDE serves as control.
The spatial domain $\Omega \subset \mathbb{R}^d$, $d=1,2,3$, is supposed to be bounded 
and Lipschitz. 
We assumed that $0 < \nu_1 \le \nu(x,t) \le \nu_2$ for almost all $(x,t) \in Q$
with positive constants $\nu_1$ and $\nu_2$.

This standard setting was already investigated in the famous book by J.L.~Lions 
\cite{LSSC2021:LS:Lions:1971}.
Since the state equation~\eqref{LS:eqn:IBVP} has a unique solution 
$y \in Y_0$, one can conclude the existence of a unique control $u \in U$ minimizing 
the quadratic cost functional $J(S(u),u)$, where $S$ is the solution operator mapping $u \in U$
to the unique solution $y \in Y_0$ of \eqref{LS:eqn:IBVP}; 
see, e.g., \cite{LSSC2021:LS:Lions:1971} and \cite{LSSC2021:LS:Troeltzsch:2010a}.
There is an huge number of publications devoted to the numerical solution 
of the optimal control problem \eqref{LS:eqn:CostFunctional}--\eqref{LS:eqn:IBVP}
with the standard $L_2(Q)$ regularization; see, e.g.,  \cite{LSSC2021:LS:Troeltzsch:2010a}.
The overwhelming  majority of the publications uses some time-stepping
or discontinuous Galerkin method
for the time discretization  in combination 
with some space-discretization method like the finite element method;  see, e.g., \cite{LSSC2021:LS:Troeltzsch:2010a}.
The unique solvability of the optimal control problem can also be established  
by showing that the optimality system 
has a unique solution.
In \cite{LSSC2021:LS:LangerSteinbachTroeltzschYang:2020b}, the Banach-Ne\u{c}as-Babu\u{s}ka 
theorem 
was applied to the optimality system to show its well-posedness.  
Furthermore, the discrete inf-sup condition, 
which
does not follow from the inf-sup condition 
in the infinite-dimensional setting, was established for 
continuous space-time  finite element discretization
on fully unstructured simplicial space-time meshes.
The discrete inf-sup condition implies stability of the discretization 
and a priori discretization error estimates. 
Distributed controls $u$ from the space $U = L_2(0,T;H^{-1}(\Omega))$ together with 
energy regularization were investigated in \cite{LSSC2021:LS:LangerSteinbachTroeltzschYang:2020c},
where one can also find a comparison of the energy regularization with the $L_2(Q)$
and the sparse regularizations.

In this paper, we make use of the maximal parabolic regularity of the reduced optimality 
system in the case of the $L_2(Q)$ regularization and under additional assumptions imposed on the 
coefficient $\nu$. Then we can derive a stabilized finite element discretization 
of the reduced optimality system in the same way as it was done for the state equation 
in our preceding papers 
\cite{LSSC2021:LangerNeumuellerSchafelner:2019a}
The properties of the finite element scheme lead to a priori discretization
error estimates that are confirmed by the numerical experiments.
%
%
\section{Space-Time Finite Element Discretization}
\label{LS:sec:SpaceTimeFEDiscretization}
Eliminating the control $u$ from the optimality system 
by means of the gradient equation $p + \varrho u = 0$,
we arrive at the reduced optimality system 
the weak form of which reads as follows:
Find the state 
$y \in Y_0$
and the adjoint state 
$p \in P_T$
such that, for $v, q \in V = L_2(0,T; H^1_0(\Omega))$, it holds
\begin{equation}
   \label{LS:eqn:OptimalitySystemWeakForm}
   \begin{array}{rcl}
      \displaystyle
      \varrho \int_Q \Big[\partial_t y \, v + \nu \,\nabla_x y \cdot \nabla_x v \Big] dQ + \int_Q p \, v \, dQ 
       & = & 0,            \\
      \displaystyle
      - \int_Q y \, q \, dQ + \int_Q  \Big[ - \partial_t p \, q + \nu \, \nabla_x p \cdot \nabla_x q \Big] dQ
       & = & \displaystyle
      - \int_Q y_d \, q \, dQ,
   \end{array}
\end{equation}
where
$
   P_T:= \{p \in L^2(0,T;H^1_0(\Omega)):\, \partial_t p \in L^2(0,T;H^{-1}(\Omega)), p = 0 \; \mbox{on} \; \Sigma_T \}.
$
The variational reduced optimality system \eqref{LS:eqn:OptimalitySystemWeakForm} is well-posed;
see \cite[Theorem 3.3]{LSSC2021:LS:LangerSteinbachTroeltzschYang:2020b}.
Moreover, we additionally assume that the coefficient $\nu(x,t)$ is of 
bounded variation in $t$ for almost all $x \in \Omega$.
Then $\partial_t u$ and $Lu := - \mathrm{div}_x(\nu\, \nabla_{x} u )$
as well as
$\partial_t p$ and $Lp := - \mathrm{div}_x(\nu\, \nabla_{x} p )$
belong to $L_2(Q)$; see \cite{LSSC2021:LS:Dier:2015a}.
This property is called maximal parabolic regularity.
In this case, the parabolic partial differential equations 
involved in the reduced optimality system  
\eqref{LS:eqn:OptimalitySystemWeakForm}
hold in $L_2(Q)$. 
Therefore, the solution of the reduced optimality system  
\eqref{LS:eqn:OptimalitySystemWeakForm} is equivalent to the 
solution of the following system of coupled forward and backward 
systems of parabolic PDEs: 
Find  $y \in Y_0 \cap H^{L,1}(Q)$ and $p \in P_T \cap H^{L,1}(Q)$
such that the coupled PDE optimality system
\begin{equation}
   \label{LS:eqn:OptimalitySystemStrongForm}
   \begin{array}{rcl}
      \displaystyle
      \varrho \Big[\partial_t y - \mbox{div}_x(\nu \nabla_x y) \Big] & = & - p\quad \mbox{in}\; L_2(Q),    \\
      \displaystyle
      - \partial_t p \, q - \mbox{div}_x(\nu \, \nabla_x p)          & = & y - y_d\quad \mbox{in}\; L_2(Q)
   \end{array}
\end{equation}
hold, where $H^{L,1}(Q) = \{v \in H^1(Q): Lv := -\mathrm{div}_x (\nu \nabla_{x} v) \in L_2(Q)\}$.
The coupled PDE optimality system \eqref{LS:eqn:OptimalitySystemStrongForm}
is now the starting point for the construction of the coercive finite element scheme.

Let $\mathcal{T}_h$ be a 
regular decomposition of the space-time cylinder $Q$ into simplicial elements, 
i.e., $\overline{Q} = \bigcup_{K\in\mathcal{T}_h} \overline{K}$, and  $K\cap K'=\emptyset$
for all $K$ and $K'$ from $ \mathcal{T}_h $ with $K \neq K'$;  
see, e.g., 
\cite{LS:Ciarlet:1978a}  for more details. 
On the basis of the triangulation $\mathcal{T}_h$, we define the space-time finite element spaces
\begin{eqnarray}
   \label{LS:eqn:Y0h}
   Y_{0h}   & = & \{ y_h\in C(\overline{Q})  : y_h(x_K(\cdot)) \in\mathbb{P}_k(\hat{K}),\,
   \forall  K \in \mathcal{T}_h,\, y_h=0\;\mbox{on}\; {\overline \Sigma} { \cap {\overline \Sigma}_0 }
   \}, \\
   \label{LS:eqn:PTh}
   P_{Th}   & = & \{ p_h\in C(\overline{Q})  : p_h(x_K(\cdot)) \in\mathbb{P}_k(\hat{K}),\,
   \forall  K \in \mathcal{T}_h,\, p_h=0\;\mbox{on}\; {\overline \Sigma} { \cap {\overline \Sigma}_T }
   \},
\end{eqnarray}
where $x_K(\cdot)$ denotes the map from the reference element $\hat{K}$ to the finite element 
$K \in \mathcal{T}_h$, 
and $\mathbb{P}_k(\hat{K})$ is the space of polynomials of the degree $k$ on the reference element $\hat{K}$.
For brevity of the presentation, we set $\nu$ to $1$.
The same derivation can be done for $\nu$ that fulfill the condition
$\mbox{div}_x(\nu \, \nabla_x w_h)|_K \in L_2(K)$ for all 
$w_h$ from $Y_{0h}$ or $P_{Th}$ and for all
$K \in \mathcal{T}_h$ (i.e., piecewise smooth)
in addition to the conditions imposed above.
Multiplying the first PDE in \eqref{LS:eqn:OptimalitySystemStrongForm} by 
$v_h + \lambda \partial_t v_h$ with $v_h \in Y_{0h}$, and the second one by
$q_h - \lambda \partial_t q_h$
with $q_h \in P_{Th}$, integrating over $K$,
integrating by parts in the elliptic parts where the scaling 
parameter
$\lambda$
does
not appear, and summing over all $K \in \mathcal{T}_h$,
we arrive at the variational consistency identity
\begin{equation}
   \label{LS:eqn:Consistency}
   a_h(y,p;v_h,q_h) = \ell_h(v_h,q_h) \quad \forall  (v_h,q_h) \in Y_{0h} \times P_{Th},
\end{equation}
with the combined bilinear and linear forms
\begin{eqnarray}
   \label{LS:eqn:a_h}
   \nonumber
   a_h(y,p;v,q)     & = &   \sum_{K\in\mathcal{T}_h}
   \int_K\Big[\varrho\bigl(\partial_{t}y\,v+\lambda\partial_{t}y\partial_{t}v + \nabla_{x}y\cdot\nabla_{x}v - \lambda \Delta_{x}y\,\partial_{t}v\bigr) \\ \nonumber
      &&\quad\quad + p(v + \lambda \partial_{t}v) - \partial_{t}p\,q+\lambda \partial_{t}p\partial_{t}q + \nabla_{x}p\cdot\nabla_{x}q\\
      &&\quad\quad +\lambda \Delta_{x}p\,\partial_{t}q - u(q - \lambda \partial_{t}q)\Big]\,\mathrm{d}K
   \quad\mbox{and}\\
   \label{LS:eqn:l_h}
   \ell_h(v,q) & = & - \sum_{K\in\mathcal{T}_h}  \int_K y_d (q - \lambda \partial_{t}q)\,\mathrm{d}K,
\end{eqnarray}
respectively. Now, the corresponding  consistent finite element scheme reads as follows: 
Find $(y_h,p_h) \in Y_{0h} \times P_{Th}$ such that
\begin{equation}
   \label{LS:eqn:FEM}
   a_h(y_h,p_h;v_h,q_h) = \ell_h(v_h,q_h) \quad \forall  (v_h,q_h) \in Y_{0h} \times P_{Th}.
\end{equation}
Subtracting \eqref{LS:eqn:FEM} from \eqref{LS:eqn:Consistency}, we immediately get the
Galerkin orthogonality relation
\begin{equation}
   \label{LS:eqn:GO}
   a_h(y-y_h,p-p_h;v_h,q_h) = 0 \quad \forall \, (v_h,q_h) \in Y_{0h} \times P_{Th},
\end{equation}
which is crucial for deriving discretization error estimates.
%
%
\section{Discretization Error Estimates}
\label{LS:sec:DiscretizationErrorEstimates}
We first show that the bilinear $a_h$ is coercive on $Y_{0h} \times P_{Th}$
with respect to
norm
\begin{eqnarray*}
   \|(v,q)\|_{h}^2 &=& \varrho\, \|v\|_{h,T}^2 + \|q\|_{h,0}^2
   =
   \varrho\,\bigl( \|v(\cdot,T)\|_{L_2(\Omega)}^2 + \|\nabla_{x}v\|_{L_2(Q)}^2 + \lambda \|\partial_t v\|_{L_2(Q)}^2
   \bigr)\\
   && 
   \hspace*{30mm}
   + \|q(\cdot,0)\|_{L_2(\Omega)}^2 + \|\nabla_{x}q\|_{L_2(Q)}^2 + \lambda \|\partial_t q\|_{L_2(Q)}^2.
\end{eqnarray*}
Indeed, for all 
$(v_h,q_h) \in Y_{0h} \times P_{Th}$,
we get the estimate 
\begin{eqnarray}
   \label{LS:eqn:Coercivity}
   \nonumber
   a_h(v_h,q_h;v_h,q_h) &=& \sum_{K\in\mathcal{T}_h}
   \int_K\Big[\varrho\bigl(\partial_{t}v_h\,v_h+\lambda |\partial_{t}v_h|^2 + |\nabla_{x}v_h|^2 - \lambda \Delta_{x}v_h\,\partial_{t}v_h\bigr)\\ \nonumber
      &&\quad\quad + q_h(v _h+ \lambda \partial_{t}v_h) - \partial_{t}q_h\,q_h+\lambda |\partial_{t}q_h|^2 + |\nabla_{x}q_h|^2\\ \nonumber
      &&\quad\quad +\lambda \Delta_{x}q_h\,\partial_{t}q_h - v_h(q_h - \lambda \partial_{t}q_h)\Big]\,\mathrm{d}K\\ 
   &\ge& \mu_c \, \|(v_h,q_h)\|_h^2,
\end{eqnarray}
with $\mu_c = 1/2$ provided that $\lambda \le c_{inv}^{-2} h^2$, where $c_{inv}$ denotes 
the constant in the inverse inequality 
$\|\mbox{div}_x(\nabla_x w_h)\|_{L_2(K)}\le c_{inv} h^{-1} \|\nabla_x w_h\|_{L_2(K)}$
that holds for all $w_h \in Y_{0h}$ or  $w_h \in P_{Th}$.
For $k=1$, the terms  $\Delta_{x}v_h$ and $\Delta_{x}q_h$ are zero,
and we do not need the inverse inequality, but $\lambda$ should be also $O(h^2)$
in order to get an optimal convergence rate estimate.
The coercivity of the bilinear form $a_h$ immediately implies 
uniqueness and existence of the finite element solution $(y_h,p_h) \in Y_{0h} \times P_{Th}$
of \eqref{LS:eqn:FEM}. In order to prove discretization error estimates, 
we need the boundedness of the bilinear form
\begin{equation}
   \label{LS:eqn:Boundedness}
   |a_h(y,p;v_h,q_h)| \le \mu_b \|(y,p)\|_{h,*}\|(v_h,q_h)\|_{h} 
   \quad \forall (v_h,q_h) \in Y_{0h} \times P_{Th},
\end{equation}
and for all $y \in Y_{0h} + Y_0 \cap H^{L,1}(Q)$ and $p \in P_{Th} + P_T \cap H^{L,1}(Q)$,
where
\begin{eqnarray*}
   \|(y,p)\|_{h,*}^2 &=& \|(y,p)\|_{h}^2 
   + \varrho  \sum_{K\in\mathcal{T}_h}\lambda\,\|\Delta_x y\|_{L_2(K)}^2 
   + [(\varrho + 1)\lambda^{-1} +\lambda]\, \|y\|_{L_2(Q)}^2\\
   && +  \sum_{K\in\mathcal{T}_h}\lambda\,\|\Delta_x p\|_{L_2(K)}^2 
   + [2\lambda^{-1} +\lambda]\, \|p\|_{L_2(Q)}^2\\
\end{eqnarray*}
Indeed, using Cauchy's inequalities and the Friedrichs inequality 
$\|w\|_{L_2(Q)} \le c_{F\Omega} \|\nabla_x w\|_{L_2(Q)}$
that holds for all $w \in Y_{0}$ or  $w \in P_{T}$,
we can easily prove \eqref{LS:eqn:Boundedness} with 
$\mu_b = (\max\{4, 1 + \lambda c_{F\Omega}^2, 3 + \varrho^{-1}, 1 + \lambda c_{F\Omega}^2 \varrho^{-1}\})^{1/2}$.
Now, \eqref{LS:eqn:GO}, \eqref{LS:eqn:Coercivity}, and \eqref{LS:eqn:Boundedness}
immediately lead to the following C\'{e}a-like estimate of the discretization 
error by some best-approximation error.
\begin{theorem}
   \label{LS:the:CEA}
   Let $y_d \in L_2(Q)$ be a given target, and let $\nu \in L_\infty(Q)$ fulfill 
   the assumptions imposed above. Furthermore, we assume that the regularization (cost)
   parameter $\varrho \in \mathbb{R}_+$ is fixed.
   Then the C\'{e}a-like estimate
   \begin{equation*}
      \|(y-y_h,p-p_h)\|_h \le \inf_{v_h \in Y_{0h}, q_h \in P_{Th}}
      \Bigl( \|(y-v_h,p-q_h)\|_h  + \frac{\mu_b}{\mu_c} \|(y-v_h,p-q_h)\|_{h,*}\Bigr)
   \end{equation*}
   holds, where $(y,p)$ and $(y_h,p_h)$ are the solutions of \eqref{LS:eqn:OptimalitySystemWeakForm}
   and \eqref{LS:eqn:FEM}, respectively.
\end{theorem}
This C\'{e}a-like estimate immediately yields convergence rate estimates of the form 
\begin{equation}
   \|(y-y_h,p-p_h)\|_h \le c(u,p) h^s 
\end{equation}
with $s = \min\{k,l\}$ provided that $y \in Y_0 \cap H^{L,1}(Q) \cap H^{l+1}(Q)$
and $p \in P_T \cap H^{L,1}(Q) \cap H^{l+1}(Q)$, where $l$ is some positive real number 
defining the regularity of the solution; see 
\cite{LSSC2021:LangerNeumuellerSchafelner:2019a}
for corresponding convergence rate estimates for the state equation only.
%
%
\section{Numerical Results}
\label{LS:sec:Numerical Result}
Let $ \{ \phi^{(j)} : j=1,\dots,N_h \} $ be a nodal finite element basis for $ Y_{0h} $, and let $ \{ \psi^{(m)} : m=1,\dots,M_h \} $ be a nodal finite element basis for $ P_{Th} $. Then we can express each finite element function $ y_h \in Y_{0h} $ and $ p_h \in P_{Th}$ via the finite element basis, i.e. $ y_h = \sum_{j=1}^{N_h} y_{j} \phi^{(j)} $ and $ p_h = \sum_{m=1}^{M_h} p_{m} \psi^{(m)} $, respectively. We insert this ansatz into \eqref{LS:eqn:FEM}, test with basis functions $ \phi^{(i)} $ and $\psi^{(n)}$, and obtain the system
\[
   \mathbf{K}_h \begin{pmatrix}
      \mathbf{y}_h \\\mathbf{p}_h
   \end{pmatrix}
   =\begin{pmatrix}
      \mathbf{0} \\
      \mathbf{f}_h
   \end{pmatrix}
\]
with 
$ \mathbf{K}_{h} = (a_h(\phi^{(j)},\psi^{(m)};\phi^{(i)},\psi^{(n)}))_{i,j=1,\dots,N_h}^{m,n=1,\dots,M_h} $,
$ \mathbf{f}_h = (\ell_h(0,\psi^{(n)}))_{n=1,\dots,M_h} $, $\mathbf{y}_h = (y_j)_{j=1,\dots,N_h}$ and $\mathbf{p}_h = (p_m)_{m=1,\dots,M_h}$. The (block)-matrix $ \mathbf{K}_h$ is non-symmetric, but positive definite due to \eqref{LS:eqn:Coercivity}. Hence the linear system is solved by means of the flexible General Minimal Residual (GMRES) method, preconditioned by a block-diagonal algebraic multigrid (AMG) method, i.e., we apply an AMG preconditioner to each of the diagonal blocks of $ \mathbf{K}_h$. Note that we only need to solve once in order to obtain a numerical solution of the space-time tracking optimal control problem \eqref{LS:eqn:CostFunctional}--\eqref{LS:eqn:IBVP}, consisting of state and adjoint state. The control can then be recovered from the gradient equation $ p + \varrho u = 0 $.

The space-time finite element method is implemented 
by means of
the \texttt{C++} library MFEM \cite{LS:mfem-library}. We use \emph{BoomerAMG}, provided by the linear solver library \emph{hypre}, to realize the preconditioner. The linear solver is stopped once the initial residual is reduced by a factor of $10^{-8}$. We are 
interested in convergence rates
with respect to
the mesh size $h$
for a fixed regularization (cost) parameter $\varrho$.
%
%
\subsection{Smooth Target}
\begin{figure}[!htb]
   \centering%
   \includegraphics[width=\linewidth]{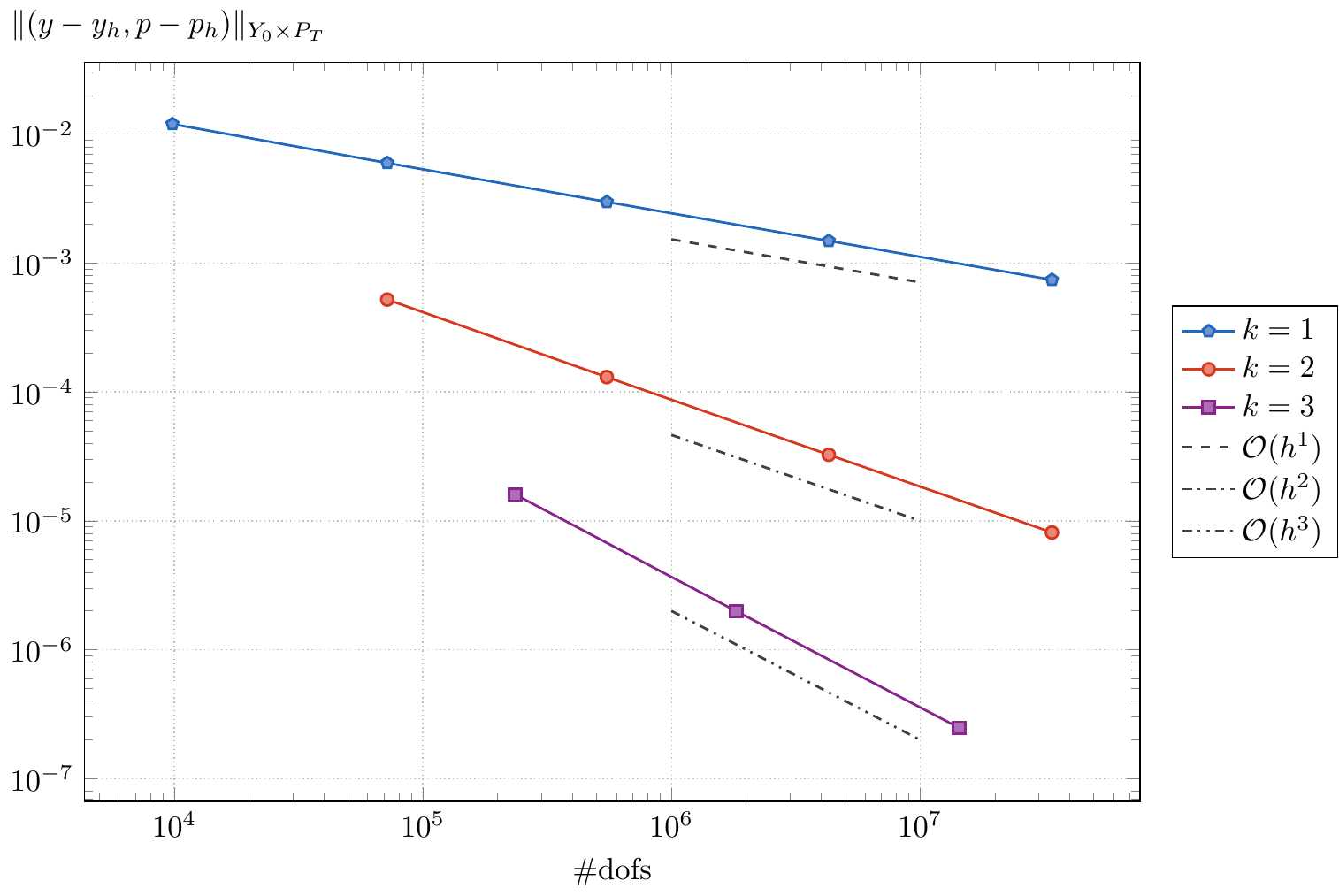}%
   \caption{Convergence rates for different polynomial degrees $ k = 1,2,3 $.}\label{LS:fig:conv}
\end{figure}

For our first example, we consider the space-time cylinder $ Q = (0,1)^3$, i.e., $d=2$, the manufactured state
\[
   y(x,t) = \sin(x_1\,\pi)\sin(x_2\,\pi) \left(a\,t^2+b\,t\right),
\]
as well as the corresponding adjoint state
\[
   p(x,t) = -\varrho \sin(x_1\,\pi)\sin(x_2\,\pi) \left(2\,\pi^2\,a\,t^2 + (2\,\pi^2\,b + 2\,a)t + b\right),
\]
with $ a=\frac{2\,\pi^2 + 1}{2\,\pi^2 + 2}\ \text{and}\ b=1. $
The desired state $ y_d $ and the optimal control $ u $ are then computed accordingly, and we fix the regularization parameter $ \varrho = 0.01 $. This problem is very smooth and devoid of any local features or singularities, hence we expect optimal convergence rates. Indeed, as we can observe in Fig.~\ref{LS:fig:conv}, the error in the $ \|(\cdot,\cdot)\|_{Y_0\times P_T} $-norm decreases with a rate of $ \mathcal{O}(h^{k}) $, where $ k $ is the polynomial degree of the finite element basis functions.
%
\subsection{Discontinuous Target}
For the second example, we consider once more the space-time cylinder $ Q = (0,1)^3$, and specify the target state
\[
   y_d(x,t) = \begin{cases}
      1, & \sqrt{(x_1 - 0.5)^2 + (x_2 - 0.5)^2 + (t-0.5)^2} \le 0.25, \\
      0, & \text{else},
   \end{cases}
\]
as an expanding and 
shrinking
circle that is nothing but a \emph{fixed} ball in the space-time
cylinder $Q$.
We use the fixed regularization parameter $ \varrho = 10^{-6} $. Here, we do not know the exact solutions for the state or the optimal control, thus we cannot consider any convergence rates for the discretization error. However, the discontinuous target state may introduce local features at the (hyper-)surface of discontinuity. Hence it might be beneficial to use adaptive mesh refinements driven by an a posteriori error indicator. 
In particular, we use the residual based indicator proposed by Steinbach and Yang \cite{LS:SteinbachYang:2018a}, applied to the residuals of the reduced optimality system \eqref{LS:eqn:OptimalitySystemStrongForm}. 
The final indicator is then the sum of the squares of both parts.

\begin{figure}[!htb]
   \centering%
   \includegraphics[width=\linewidth]{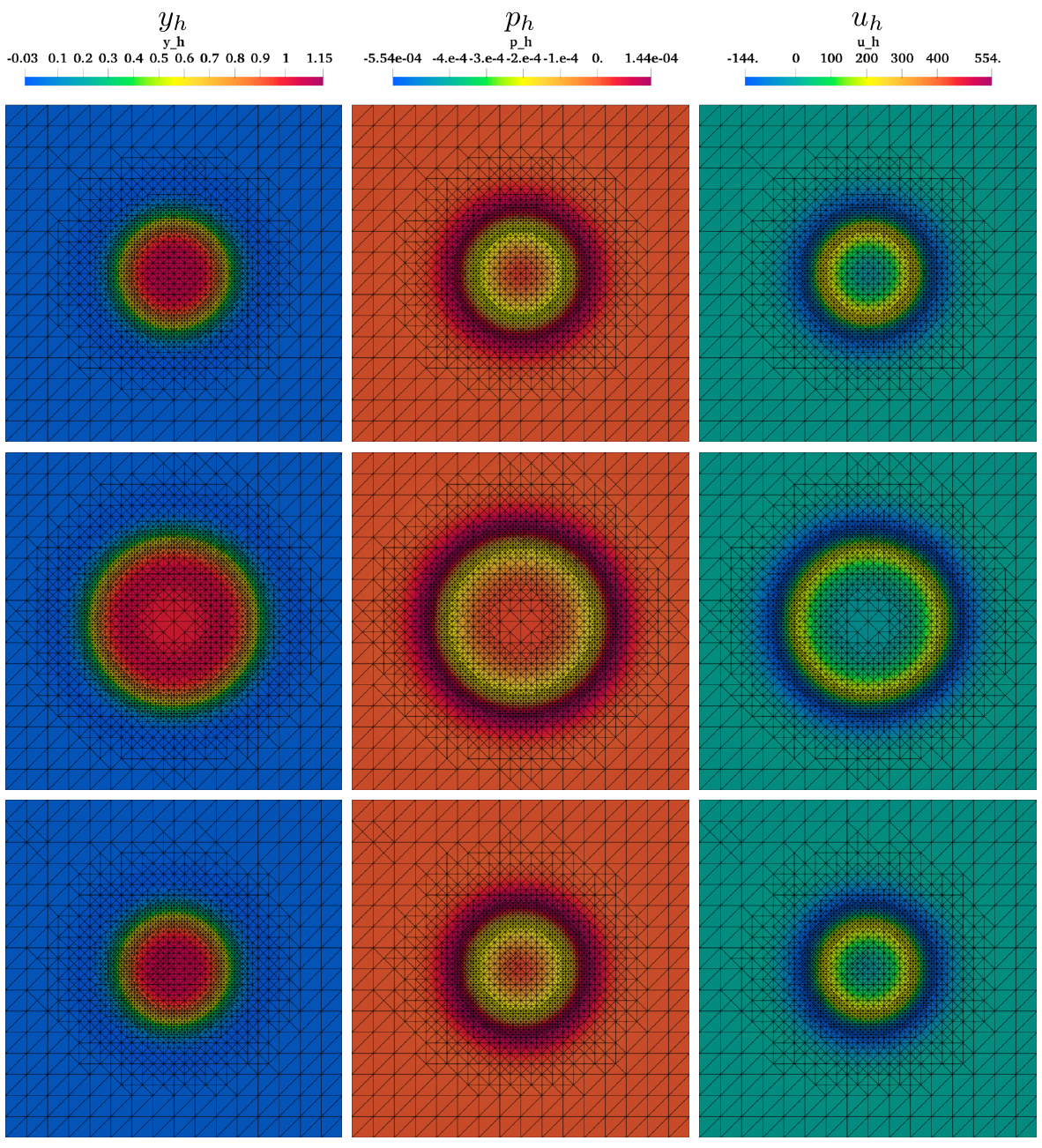}%
   \caption{Finite element solutions, with $ J(y_h,u_h) = 3.5095\times10^{-3} $, plotted over the space-time mesh, obtained after 20 adaptive refinements, and cut at $t = 0.3125$ (upper row), $t=0.5$ (middle row), and $t=0.6875$ (lower row).}
   \label{LS:fig:meshes}
\end{figure}

In Fig.~\ref{LS:fig:meshes}, we present the finite element functions $y_h$, $p_h$, and $u_h$, plotted over cuts of the space-time mesh $ \mathcal{T}_h $ at different times $t$. We can observe that the mesh refinements are mostly concentrated in annuli centered at $ (0.5, 0.5) $, e.g. for $ t = 0.5$, the outer and inner radii are $ \sim\frac{7}{36}\pm \frac{1}{36} $, respectively; see Fig.~\ref{LS:fig:meshes} (middle row). 
%
%
\section{Conclusions}
\label{LS:sec:Conclusions}
We proposed a stable, fully unstructured, space-time simplicial finite element discretization of the reduced 
optimality system of the standard space-time tracking parabolic optimal control problem 
with $L_2$-regularization.
We derived a priori discretization error estimates. We presented numerical results 
for two benchmarks. We observed optimal rates for the example with smooth solutions 
as predicted by the a priori estimates. In the case of a discontinuous target,
we use full space-time adaptivity. In order to get the full space-time solution $(y_h,p_h,u_h)$,
one has to solve only {\it one} system of algebraic equations. 
In this paper, we used flexible GMRES preconditioned by AMG. 
\bibliographystyle{acm}
\bibliography{ms}
\end{document}